# A CLASS OF RECURSIVE SETS

Florentin Smarandache
University of New Mexico
200 College Road
Gallup, NM 87301, USA
E-mail: smarand@unm.edu

In this article one builds a class of recursive sets, one establishes properties of these sets and one proposes applications. This article widens some results of [1].

**1) Definitions, properties.**

One calls recursive sets the sets of elements which are built in a recursive manner: let $T$ be a set of elements and $f_i$ for $i$ between 1 and $s$, of operations $n_i$, such that $f_i : T^{n_i} \to T$. Let's build by recurrence the set $M$ included in $T$ and such that:

**(Def. 1)** $1^o$) certain elements $a_1, ..., a_n$ of $T$, belong to $M$.

$2^o$) if $(\alpha_{i_1}, ..., \alpha_{i_{n_i}})$ belong to $M$, then $f_i(\alpha_{i_1}, ..., \alpha_{i_{n_i}})$ belong to $M$ for all $i \in \{1, 2, ..., s\}$.

$3^o$) each element of $M$ is obtained by applying a number finite of times the rules $1^o$ or $2^o$.

We will prove several proprieties of these sets $M$, which will result from the manner in which they were defined. The set $M$ is the representative of a class of recursive sets because in the rules $1^o$ and $2^o$, by particularizing the elements $a_1, ..., a_n$ respectively $f_1, ..., f_s$ one obtains different sets.

**Remark 1:** To obtain an element of $M$, it is necessary to apply initially the rule 1.

**(Def. 2)** The elements of $M$ are called elements $M$-recursive.

**(Def. 3)** One calls order of an element $a$ of $M$ the smallest natural $p \geq 1$ which has the propriety that $a$ is obtained by applying $p$ times the rule $1^o$ or $2^o$.

One notes $M_p$ the set which contains all the elements of order $p$ of $M$. It is obvious that $M_1 = \{a_1, ..., a_n\}$.

$$M_2 = \bigcup_{i=1}^{s} \left\{ \bigcup_{(\alpha_{i_1}, ..., \alpha_{i_{n_i}}) \in M_1^{n_i}} f_i(\alpha_{i_1}, ..., \alpha_{i_{n_i}}) \right\} \setminus M_1.$$

One withdraws $M_1$ because it is possible that $f_j(a_{j_1}, ..., a_{j_{n_j}}) = a_i$ which belongs to $M_1$, and thus does not belong to $M_2$.

One proves that for $k \geq 1$ one has:

$$M_{k+1} = \bigcup_{i=1}^{s} \left\{ \bigcup_{(\alpha_{i_1},...,\alpha_{i_{n_i}}) \in \prod_k^{(i)}} f_i(\alpha_{i_1},...,\alpha_{i_{n_i}}) \right\} \setminus \bigcup_{h=1}^{k} M_h$$

where each

$$\prod_k^{(i)} = \{(\alpha_{i_1},...,\alpha_{i_{n_i}}) / \alpha_{i_j} \in M_{q_j} \quad j \in \{1,2,...,n_i\}; \; 1 \leq q_j \leq k \text{ and at least an}$$

element $a_{i_{j_o}} \in M_k, 1 \leq j_o \leq n_i \}$.

The sets $M_p$, $p \in \mathbb{N}^*$, form a partition of the set $M$.

**Theorem 1:**

$$M = \bigcup_{p \in \mathbb{N}^*} M_p, \text{ where } \mathbb{N}^* = \{1, 2, 3,...\}.$$

*Proof:*
From the rule 1º it results that $M_1 \subseteq M$.

One supposes that this propriety is true for values which are less than $p$. It results that $M_p \subseteq M$, because $M_p$ is obtained by applying the rule 2º to the elements of $\bigcup_{i=1}^{p-1} M_i$.

Thus $\bigcup_{p \in \mathbb{N}^*} M_p \subseteq M$. Reciprocally, one has the inclusion in the contrary sense in accordance with the rule 3º.

**Theorem 2:** The set $M$ is the smallest set, which has the properties 1º and 2º.
*Proof:*
Let $R$ be the smallest set having properties 1º and 2º. One will prove that this set is unique.

Let's suppose that there exists another set $R'$ having properties 1º and 2º, which is the smallest. Because $R$ is the smallest set having these proprieties, and because $R'$ has these properties also, it results that $R \subseteq R'$; of an analogue manner, we have $R' \subseteq R$: therefore $R = R'$.

It is evident that $M' \subseteq R$. One supposes that $M_i \subseteq R$ for $1 \leq i < p$. Then (rule 3º), and taking in consideration the fact that each element of $M_p$ is obtained by applying rule 2º to certain elements of $M_i$, $1 \leq i < p$, it results that $M_p \subseteq R$. Therefore $\bigcup_p M_p \subseteq R$ $(p \in \mathbb{N}^*)$, thus $M \subseteq R$. And because $R$ is unique, $M = R$.

**Remark 2.** The theorem 2 replaces the rule 3º of the recursive definition of the set $M$ by: " $M$ is the smallest set that satisfies proprieties 1º and 2º".

**Theorem 3:** $M$ is the intersection of all the sets of $T$ which satisfy conditions 1º and 2º.
*Proof:*

Let $T_{12}$ be the family of all sets of $T$ satisfying the conditions 1° and 2°. We note $I = \bigcap_{A \in T_{12}} A$.

$I$ has the properties 1° and 2° because:
1) For all $i \in \{1,2,...,n\}$, $a_i \in I$, because $a_i \in A$ for all $A$ of $T_{12}$.
2) If $\alpha_{i_1},...,\alpha_{i_{n_i}} \in I$, it results that $\alpha_{i_1},...,\alpha_{i_{n_i}}$ belong to $A$ that is $A$ of $T_{12}$. Therefore,
$\forall i \in \{1,2,...,s\}$, $f_i(\alpha_{i_1},...,\alpha_{i_{n_i}}) \in A$ which is $A$ of $T_{12}$, therefore $f_i(\alpha_{i_1},...,\alpha_{i_{n_i}}) \in I$ for all $i$ from $\{1,2,...,s\}$.

From theorem 2 it results that $M \subseteq I$.

Because $M$ satisfies the conditions 1° and 2°, it results that $M \in T_{12}$, from which $I \subseteq M$. Therefore $M = I$

**(Def. 4)** A set $A \subseteq I$ is called closed for the operation $f_{i_0}$ if and only if for all $\alpha_{i_0 1},...,\alpha_{i_0 n_{i_0}}$ of $A$, one has $f_{i_0}(\alpha_{i_0 1},...,\alpha_{i_0 n_{i_0}})$ belong to $A$.

**(Def. 5)** A set $A \subseteq T$ is called $M$-recursively closed if and only if:
1) $\{a_1,...,a_n\} \subseteq A$.
2) $A$ is closed in respect to operations $f_1,...,f_s$.

With these definitions, the precedent theorems become:

**Theorem 2':** The set $M$ is the smallest $M$-recursively closed set.

**Theorem 3':** $M$ is the intersection of all $M$-recursively closed sets.

**(Def. 6)** The system of elements $\langle \alpha_1,...,\alpha_m \rangle$, $m \geq 1$ and $\alpha_i \in T$ for $i \in \{1,2,...,m\}$, constitute a $M$-recursive description for the element $\alpha$, if $\alpha_m = \alpha$ and that each $\alpha_i$ ($i \in \{1,2,...,m\}$) satisfies at least one of the proprieties:
1) $\alpha_i \in \{a_1,...,a_n\}$.
2) $\alpha_i$ is obtained starting with the elements which precede it in the system by applying the functions $f_j$, $1 \leq j \leq s$ defined by property 2° of (Def. 1).

**(Def. 7)** The number $m$ of this system is called the length of the $M$-recursive description for the element $\alpha$.

**Remark 3:** If the element $\alpha$ admits a $M$-recursive description, then it admits an infinity of such descriptions.

Indeed, if $\langle \alpha_1,...,\alpha_m \rangle$ is a $M$-recursive description of $\alpha$ then $\langle \underbrace{a_1,...,a_1}_{h \text{ times}}, \alpha_1,...,\alpha_m \rangle$ is also a $M$-recursive description for $\alpha$, $h$ being able to take all values from $\mathbb{N}$.

**Theorem 4:** The set $M$ is identical with the set of all elements of $T$ which admit a $M$-recursive description.

*Proof:* Let $D$ be the set of all elements, which admit a $M$-recursive description. We will prove by recurrence that $M_p \subseteq D$ for all $p$ of $\mathbb{N}^*$.

For $p = 1$ we have: $M_1 = \{a_1,...,a_n\}$, and the $a_j$, $1 \leq j \leq n$, having as $M$-recursive description: $\langle a_j \rangle$. Thus $M_1 \subseteq D$. Let's suppose that the property is true for the values smaller than $p$. $M_p$ is obtained by applying the rule 2° to the elements of $\bigcup_{i=1}^{p-1} M_i$; $\alpha \in M_p$ implies that $\alpha \in f_j(\alpha_{i_1},...,\alpha_{i_{n_i}})$ and $\alpha_{i_j} \in M_{h_j}$ for $h_j < p$ and $1 \leq j \leq n_i$.

But $a_{i_j}$, $1 \leq j \leq n_i$, admits $M$-recursive descriptions according to the hypothesis of recurrence, let's have $\langle \beta_{j1},...,\beta_{js_j} \rangle$. Then $\langle \beta_{11},...,\beta_{1s_1},\beta_{21},...,\beta_{2s_2},...,\beta_{n_i 1},...,\beta_{n_i s_{n_i}},\alpha \rangle$ constitute a $M$-recursive description for the element $\alpha$. Therefore if $\alpha$ belongs to $D$, then $M_p \subseteq D$ which is $M = \bigcup_{p \in \mathbb{N}^*} M_p \subseteq D$.

Reciprocally, let $x$ belong to $D$. It admits a $M$-recursive description $\langle b_1,...,b_t \rangle$ with $b_t = x$. It results by recurrence by the length of the $M$-recursive description of the element $x$, that $x \in M$. For $t = 1$ we have $\langle b_1 \rangle$, $b_1 = x$ and $b_1 \in \{a_1,...,a_n\} \subseteq M$. One supposes that all elements $y$ of $D$ which admit a $M$-recursive description of a length inferior to $t$ belong to $M$. Let $x \in D$ be described by a system of length $t$: $\langle b_1,...,b_t \rangle$, $b_t = x$. Then $x \in \{a_1,...,a_n\} \subseteq M$, where $x$ is obtained by applying the rule 2° to the elements which precede it in the system: $b_1,...,b_{t-1}$. But these elements admit the $M$-recursive descriptions of length which is smaller that $t$: $\langle b_1 \rangle, \langle b_1, b_2 \rangle,...,\langle b_1,...,b_{t-1} \rangle$. According to the hypothesis of the recurrence, $b_1,...,b_{t-1}$ belong to $M$. Therefore $b_t$ belongs also to $M$. It results that $M \equiv D$.

**Theorem 5:** Let $b_1,...,b_q$ be elements of T, which are obtained from the elements $a_1,...,a_n$ by applying a finite number of times the operations $f_1, f_2,...$, or $f_s$. Then $M$ can be defined recursively in the following mode:
1) Certain elements $a_1,...,a_n, b_1,...,b_q$ of $T$ belong to $M$.
2) $M$ is closed for the applications $f_i$, with $i \in \{1, 2,..., s\}$.
3) Each element of $M$ is obtained by applying a finite number of times the rules (1) or (2) which precede.

*Proof:* evident. Because $b_1,...,b_q$ belong to $T$, and are obtained starting with the elements $a_1,...,a_n$ of $M$ by applying a finite number of times the operations $f_i$, it results that $b_1,...,b_q$ belong to $M$.

**Theorem 6:** Let's have $g_j$, $1 \le j \le r$, of the operations $n_j$, where $g_j : T^{n_j} \to T$ such that $M$ to be closed in rapport to these operations. Then $M$ can be recursively defined in the following manner:
1) Certain elements $a_1,...,a_n$ de $T$ belong to $M$.
2) $M$ is closed for the operations $f_i$, $i \in \{1,2,...,s\}$ and $g_j$, $j \in \{1,2,...,r\}$.
3) Each element of $M$ is obtained by applying a finite number of times the precedent rules.

Proof is simple: Because $M$ is closed for the operations $g_j$ (with $j \in \{1,2,...,r\}$), one has, that for any $\alpha_{j1},...,\alpha_{jn_j}$ from $M$, $g_j(\alpha_{j1},...,\alpha_{jn_j}) \in M$ for all $j \in \{1,2,...,r\}$.

From the theorems 5 and 6 it results:

**Theorem 7**: The set M can be recursively defined in the following manner:
1) Certain elements $a_1,...,a_n,b_1,...,b_q$ of $T$ belong to $M$.
2) $M$ is closed for the operations $f_i$ ($i \in \{1,2,...,s\}$) and for the operations $g_j$ ($j \in \{1,2,...,r\}$) previously defined.
3) Each element of M is defined by applying a finite number of times the previous 2 rules.

**(Def. 8)** The operation $f_i$ conserves the property $P$ iff for any elements $\alpha_{i1},...,\alpha_{in_i}$ having the property $P$, $f_i(\alpha_{i1},...,\alpha_{in_i})$ has the property $P$.

**Theorem 8**: If $a_1,...,a_n$ have the property $P$, and if the functions $f_1,...,f_s$ preserve this property, then all elements of $M$ have the property $P$.

*Poof*:

$M = \bigcup\limits_{p \in \mathbb{N}^*} M_p$. The elements of $M_1$ have the property $P$.

Let's suppose that the elements of $M_i$ for $i < p$ have the property $P$. Then the elements of $M_p$ also have this property because $M_p$ is obtained by applying the operations $f_1, f_2,..., f_s$ to the elements of: $\bigcup\limits_{i=1} M_i$, elements which have the property $P$.

Therefore, for any $p$ of $\mathbb{N}$, the elements of $M_p$ have the property $P$.

Thus all elements of $M$ have it.

**Corollary 1**: Let's have the property $P$: "$x$ can be represented in the form $F(x)$".

If $a_1,...,a_n$ can be represented in the form $F(a_1),...,$ respectively $F(a_n)$, and if $f_1,...,f_s$ maintains the property $P$, then all elements $\alpha$ of $M$ can be represented in the form $F(\alpha)$.

Remark. One can find more other equivalent def. of $M$.

## 2) APPLICATIONS, EXAMPLES.

In applications, certain general notions like: $M$ - recursive element, $M$ -recursive description, $M$ - recursive closed set will be replaced by the attributes which characterize the set $M$. For example in the theory of recursive functions, one finds notions like: recursive primitive functions, primitive recursive description, primitively recursive closed sets. In this case "$M$" has been replaced by the attribute "primitive" which characterizes this class of functions, but it can be replaced by the attributes "general", "partial".

By particularizing the rules 1° and 2° of the def. 1, one obtains several interesting sets:

**Example 1:** (see [2], pp. 120-122, problem 7.97).

**Example 2:** The set of terms of a sequence defined by a recurring relation constitutes a recursive set.

Let's consider the sequence: $a_{n+k} = f(a_n, a_{n+1}, ..., a_{n+k-1})$ for all $n$ of $\mathbb{N}^*$, with $a_i = a_i^0$, $1 \leq i \leq k$. One will recursively construct the set $A = \{a_m\}_{m \in \mathbb{N}^*}$ and one will define in the same time the position of an element in the set $A$:

1°) $a_1^0, ..., a_k^0$ belong to $A$, and each $a_i^0$ ($1 \leq i \leq k$) occupies the position $i$ in the set $A$;

2°) if $a_n, a_{n+1}, ..., a_{n+k-1}$ belong to $A$, and each $a_j$ for $n \leq j \leq n+k-1$ occupies the position $j$ in the set $A$, then $f(a_n, a_{n+1}, ..., a_{n+k-1})$ belongs to $A$ and occupies the position $n+k$ in the set $A$.

3°) each element of $B$ is obtained by applying a finite number of times the rules 1° or 2°.

**Example 3:** Let $G = \{e, a^1, a^2, ..., a^p\}$ be a cyclic group generated by the element $a$. Then $(G, \square)$ can be recursively defined in the following manner:

1°) $a$ belongs to $G$.

2°) if $b$ and $c$ belong to $G$ then $b \square c$ belongs to $G$.

3°) each element of $G$ is obtained by applying a finite number of times the rules 1 or 2.

**Example 4:** Each finite set $ML = \{x_1, x_2, ..., x_n\}$ can be recursively defined (with $ML \subseteq T$):

1°) The elements $x_1, x_2, ..., x_n$ of $T$ belong to $ML$.

2°) If $a$ belongs to $ML$, then $f(a)$ belongs to $ML$, where $f: T \to T$ such that $f(x) = x$;

3°) Each element of $ML$ is obtained by applying a finite number of times the rules 1° or 2°.

**Example 5:** Let $L$ be a vectorial space on the commutative corps $K$ and $\{x_1, ..., x_m\}$ be a base of $L$. Then $L$, can be recursively defined in the following manner:

1°) $x_1, ..., x_m$ belong to $L$;

2°) if $x, y$ belong to $L$ and if $a$ belongs to $K$, then $x \perp y$ $y$ belong to $L$ and $a * x$ belongs to $L$;

3°) each element of $L$ is recursively obtained by applying a finite number of times the rules 1° or 2°.

(The operators $\perp$ and $*$ are respectively the internal and external operators of the vectorial space $L$).

**Example 6:** Let $X$ be an $A$-module, and $M \subset X$ ($M \neq \emptyset$), with $M = \{x_i\}_{i \in I}$. The sub-module generated by $M$ is:

$$\langle M \rangle = \{x \in X \,/\, x = a_1 x_1 + \ldots + a_n x_n, \ a_i \in A, \ x_i \in M, \ i \in \{1, \ldots, n\}\}$$

can be recursively defined in the following way:

1°) for all $i$ of $\{1,2,\ldots,n\}$, $\{1,2,\ldots,n\} \Box x_i \in \langle M \rangle$;

2°) if $x$ and $y$ belong to $\langle M \rangle$ and $a$ belongs to $A$, then $x+y$ belongs to $\langle M \rangle$, and $ax$ also;

3°) each element of $\langle M \rangle$ is obtained by applying a finite number of times the rules 1° or 2°.

In accordance to the paragraph 1 of this article, $\langle M \rangle$ is the smallest sub-set of X that verifies the conditions 1° and 2°, that is $\langle M \rangle$ is the smallest sub-module of X that includes $M$. $\langle M \rangle$ is also the intersection of all the subsets of $X$ that verify the conditions 1° and 2°, that is $\langle M \rangle$ is the intersection of all sub-modules of $X$ that contain $M$. One also directly refines some classic results from algebra.

One can also talk about sub-groups or ideal generated by a set: one also obtains some important applications in algebra.

**Example 7:** One also obtains like an application the theory of formal languages, because, like it was mentioned, each regular language (linear at right) is a regular set and reciprocally. But a regular set on an alphabet $\Sigma = \{a_1, \ldots, a_n\}$ can be recursively defined in the following way:

1°) $\emptyset, \{\varepsilon\}, \{a_1\}, \ldots, \{a_n\}$ belong to $R$.

2°) if $P$ and $Q$ belong to $R$, then $P \cup Q$, $PQ$, and $P^*$ belong to $R$, with

$$P \cup Q = \{x \,/\, x \in P \text{ or } x \in Q\}; \quad PQ = \{xy \,/\, x \in P \text{ and } y \in Q\}, \quad \text{and} \quad P^* = \bigcup_{n=0}^{\infty} P^n \quad \text{with}$$

$P^n = \underbrace{P \cdot P \cdots P}_{n \text{ times}}$ and, by convention, $P^0 = \{\varepsilon\}$.

3°) Nothing else belongs to $R$ other that those which are obtained by using 1° or 2°.

From which many properties of this class of languages with applications to the programming languages will result.